\documentclass[12pt]{article}
\usepackage{latexsym,amsmath,amsthm,verbatim,ifthen,amssymb}
\usepackage[all]{xy}

\newtheorem{theorem}{Theorem}
\newtheorem{corollary}[theorem]{Corollary}
\newtheorem{lemma}[theorem]{Lemma}
\newtheorem{proposition}[theorem]{Proposition}

\newtheorem{remark}[theorem]{Remark}

\newcommand{\locp}{{(p)}}

\newcommand{\bz}{\mathbb Z}

\newcommand{\bq}{\mathbb Q}

\newcommand{\aut}{\text{aut}\,}

\newcommand{\eva}{\mathcal E}

\newcommand{\es}{\mathcal E_{ \sharp }}
\newcommand{\esp}{\mathcal E_{ \sharp p }}

\newcommand{\nil}{\text{nil}\,}

\newcommand{\calm}{{\cal M}}
\newcommand{\calh}{{\cal H}}
\newcommand{\ecalm}{{\mathcal E}(X;\calm)}

\begin{document}

\title{Nilpotency of self homotopy equivalences with coefficients}
\author{Maxence Cuvilliez, Aniceto Murillo\footnote{The second author
acknowledges the kind hospitality and support from the {\em Fields
Institute of the University of Toronto}
where part of the research in this paper was carried on.} and
Antonio Viruel\footnote{The second and third author have been
partially supported by the Ministerio de Educaci\'on y Ciencia
grant MTM2007-60016, and Junta de Andaluc{\'\i}a grants FQM-213
and P07-FQM-2863.\hfill\break 2000 M.S.C.: 55P10,
55P60.\hfill\break Keywords: Self homotopy equivalence.}\\
Departamento de \'Algebra, Geometr\'\i a y Topolog\'\i a,\\
Universidad de M\'alaga,\\ Ap. 59, 29080 M\'alaga, SPAIN.}

\maketitle

\begin{abstract} In this paper we study the nilpotency of certain
groups of self homotopy equivalences. Our main goal is
to extend, to localized homotopy groups and/or homotopy groups
with coefficients, the general principle of Dror and Zabrodsky by
which a group of self homotopy equivalences of a finite space
which acts nilpotently on the homotopy groups is itself nilpotent.
\end{abstract}

\section{Introduction}

Given a pointed space $X$, denote by $\eva(X)$ the group of
(based) self homotopy equivalences, i.e., the group of
automorphisms of $X$ in the pointed homotopy category. From now on
we shall consider connected complexes of finite type $X$ which are
either finite or with finitely many non trivial homotopy groups.
We denote by $\dim X=N$ its topological or homotopical dimension.
Unless explicitly stated otherwise, all  spaces will be of this
kind.

Although the computation of  $\eva(X)$ is known to be a hard task,
there are two classical and key results that impose to this group
important structural constraints:

On one hand, a theorem of Sullivan  \cite[Theorem 10.3]{su} and
Wilkerson   \cite[Theorem B]{wi} states that $\eva(X)$ is finitely
presented. This was originally proved for simply connected spaces
and later on generalized to virtually nilpotent spaces by Dror,
Dwyer and Kan \cite[Theorem 1.1]{ddk}. The main step in the proof
is to show that $\eva(X_\bq)$ is an algebraic group and that
$\eva(X)$ is commensurable with an arithmetic subgroup of
$\eva(X_\bq)$. As a consequence, it can be shown that there exists
a finite bound for the finite orders of elements of $\eva(X)$.

On the other hand we have the following theorem  due to Dror and
Zabrodsky:

\begin{theorem}\label{principal}{\em \cite[Theorem B]{d-z}} Let $G$ be a subgroup of
$\eva(X)$ which acts nilpotently on $\pi_{\le N}(X)$. Then $G$ is
itself nilpotent. In particular, $\es^m(X)$ is nilpotent.
\end{theorem}

Recall that, for $0\le m\le\infty$, $\es^m(X)$ is the
distinguished subgroup   of $\eva(X)$ formed by those classes
inducing the identity on the homotopy groups up to $m$. In other
words, $$ \es^m(X)=\text{ker}\bigl(\eva(X)\longrightarrow
\Pi_{i\le m}\text{aut}\,\pi_iX\bigr). $$  If $\dim X=N$ we shall
denote  $\es^N(X)$ simply by  $\es(X)$.

  Here, we present a slightly different proof for this well known result
in which we use a broader study of self homotopy equivalences in
the homotopy category ${\cal L}^*$ of (based) spaces with local
coefficients. Recall (see \cite[Chap.VI]{white} for instance) that
objects in this category are pairs $(X,\calm)$ in which $X$ is a
(based) topological space and $\calm=\{M_x\}_{x\in X}$ is a local
coefficient system in $X$. On the other hand, a morphism
$(f,\Theta)\colon (X,\calm)\to(Y,\calh)$ is a pair formed by a
based map $f\colon X\to Y$ and a morphism $\Theta\colon
f^*\calh\to\calm$  of local coefficient. By $f^*\calh$ we denote,
as usual, the local coefficient system on $X$ induced by $f$,
i.e., $(f^*\calh)_x=H_{f(x)}$. For each $x\in X$ we shall denote
by $\Theta_x\colon H_{f(x)}\to M_x$ the corresponding group
morphism at $x$. After considering the appropriate homotopy
notion, one obtains the homotopy category   ${\cal L}^*$. The
group of self homotopy equivalences of an object $(X,\calm)\in
{\cal L}^*$ shall be denoted by   $\ecalm$. Then, we prove:

\begin{theorem}\label{local} Let $X$ be a finite Postnikov piece
and let  $G\subset \ecalm$ be a subgroup which acts nilpotently
on both $\pi_*(X)$ and  $\calm$. Then, $G$ acts nilpotently on $H^*(X;\calm)$.
\end{theorem}

At the sight of Theorem \ref{principal}, and taking into account
the bound of finite orders of elements of $\eva(X)$ plus the
existence of a ``fracture lemma" for this group (see \cite[Theorem
8.2]{Rutter}), it has been of interest to study whether  $\es(X)$
satisfies the same structural restrictions when taking
$p$-localization, $p$-completion or considering $\esp(X)$. This
denotes the subgroup of $\eva(X)$ formed by those classes which
induce the identity on the homotopy groups of $X$ with
coefficients on $\bz/p$, up to the dimension of $X$, i.e., $$
\esp(X)=\text{ker}\bigl(\eva(X)\longrightarrow \Pi_{\le
N}\text{aut}\,\pi_i(X;\bz/p)\bigr). $$ As examples of this, we
mention two interesting results for a given nilpotent space $X$:
Maruyama proved \cite[Theorem 0.1]{ma} that
$\es(X)_{(p)}=\es^N(X_{(p)})$ while, on the other hand,  M{\o}ller
showed \cite[Theorem 4.3]{mo} that
$$\es(X_{\bz_p})=Ext(\bz/p^\infty,\es(X))=
\es(X_p{\displaystyle\hat{}}). $$  Here and henceforth,
$(-)_{\bz_p}$ denotes $H_*(-;\bz/p)$-localization while
$(-)_{(p)}$ and $(-)_p{\displaystyle\hat{}}$ are the classical
localization and completion on the prime $p$.

In this paper we plan to continue this investigation extending
Theorem \ref{principal} above, considering a subgroup of $\eva(X)$
which acts nilpotently in the homotopy groups of $X$ localized,
completed or with coefficients in $\bz/p$. Concerning this purpose
we prove:

\begin{theorem}\label{dos} Assume that $\pi_1(X)$ is a nilpotent group and  let $G$ be a subgroup of $\eva(X)$ which acts
nilpotently on $\pi_{\le N}(X)_{(p)}$, for $p$ any prime number
and $0$. If the nilpotency orders of all these actions are bounded
by a fixed integer, then G is nilpotent.
\end{theorem}
\begin{remark} {\em Observe that in the theorem above the condition of $\pi_1(X)$ being nilpotent is essential. Otherwise,
choose any finite simple group $G$ which is known to be
generically trivial, i.e.,  $G_{(p)}=\{1\}$ for $p$ any prime
number or zero. On the other hand observe that the map $G\to \aut
G$ given by inner automorphisms is a monomorphism. Indeed, its
kernel is the center of $G$ which is trivial since $G$ is simple.
This inclusion renders the non nilpotent group $G$ as a subgroup
of $\eva\bigl(K(G,1)\bigr)$ which acts nilpotently in the
localized homotopy group.

In a similar way, we may even produce an example of a solvable,
non nilpotent group, of homotopy equivalences acting nilpotently
in the localized homotopy groups of the space. Consider the
symmetric group $\Sigma_3$ and observe that
${\Sigma_3}_{(2)}=\bz/2$ while ${\Sigma_3}_{(p)}=1$ for $p\not=2$.
Again, $\Sigma_3\subset\eva\bigl(K(\Sigma_3,1)\bigr)$ is a
solvable non nilpotent group acting nilpotently on any
${\Sigma_3}_{(p)}$.}
\end{remark}

A more subtle and slightly  different situation is given when
considering nilpotent actions of subgroups of self homotopy
equivalences on the Frattini factor of the homotopy groups. Recall
that given a group $G$, the Frattini subgroup $\Phi(G)$ is the
intersection of all maximal proper subgroups of $G$. The quotient
$G/\Phi(G)$ is called the Frattini factor.

\begin{theorem}\label{uno} Assume that
$\pi_{\le N}(X)$ is a finite nilpotent group and let $G$ be a
subgroup of $\eva(X)$ which acts nilpotently on $\pi_{\le
N}(X)/\Phi\bigl(\pi_{\le N}(X)\bigr)$. Then, G is nilpotent.
\end{theorem}

In particular, taking into account that for an abelian $p$-group
$G$ its Frattini factor is precisely $G\otimes\bz/p$, we obtain
the following:

\begin{corollary}\label{nuevo} Assume that
$\pi_{\le N}(X)$ is a finite abelian group and let $G$ be a
subgroup of $\eva(X)$ which acts nilpotently on $\pi_{\le
N}(X)\otimes\bz/p$ for any prime $p$. Then, G is
nilpotent.\hfill$\square$
\end{corollary}

Notice that, by the Universal Coefficients Theorem for homotopy,
$\pi_*X\otimes\bz/p=\text{Ext}(\bz/p,\pi_*X)$ is a subgroup of
$\pi_*(X;\bz/p)$. Hence, as an immediate consequence of Corollary
\ref{nuevo} above we get:

\begin{corollary}\label{coeficientes} Assume that
$\pi_{\le N}(X)$ is a finite abelian group and let $G$ be a
subgroup of $\eva(X)$ which acts nilpotently on $\pi_{\le
N}(X;\bz/p)$ for  any prime $p$. Then, G is
nilpotent.\hfill$\square$
\end{corollary}

Having studied the nilpotency of a general subgroup of $\eva(X)$,
we  now focus on the group $\esp(X)$ and give necessary conditions
for it to be nilpotent.

\begin{theorem}\label{tres}  Let $X$ be a space for which $\pi_{\le N}(X)$ is a finite abelian $p$-group. Then
$\esp(X)$ is nilpotent and $\esp(X)/\es(X)$ is a finite $p$-group.

\end{theorem}

\begin{theorem}\label{cuatro}
Let $X$ be a space for which $\pi_{\le N}(X)$ is a finitely
generated abelian group. Then $\cap_{p\,\,\text{prime}}\esp(X)$ is
nilpotent.
\end{theorem}

\begin{remark}\label{importante}  {\em Observe that in general
$\esp(X)$ is bigger than $\es(X)$. For instance consider
$X=K(\bz/{p^r},n)$, $r,n\ge 2$. Obviously $\es(X)=\{1\}$, while
the automorphism $\rho$ of $\bz/{p^r}$ given by
$\rho(1)=p^{r-1}+1$ induces a non trivial element of $\esp(X)$.
Indeed, by the Universal Coefficients Theorem for homotopy, $$
\pi_*(X,\bz/p)=\pi_n(X,\bz/p)\oplus \pi_{n-1}(X,\bz/p), $$  in
which\hfill\break

\medskip
\noindent$\pi_n(X,\bz/p)=\hom(\bz/p,\bz/{p^r})$ {and}
$\pi_{n-1}(X,\bz/p)=\text{Ext}(\bz/p,\bz/{p^r})=\bz/p.$

\medskip
\noindent Trivially $\rho$ induces the identity on both. Note that
this example also shows that even
$\cap_{p\,\,\text{prime}}\esp(X)$  can be bigger than $\es(X)$.}
\end{remark}

The paper is organized as follows: in the next section we collect
the results we shall need from group theory and from which
Theorems \ref{dos} and \ref{uno} are immediately deduced. Theorem
\ref{local} and \ref{principal} are proved in section \S3.
Finally, in  section \S4 we establish Theorems \ref{tres} and
\ref{cuatro}.

\section{From group theory}

We begin by recalling some basic facts. If $G$ is a group acting
on another group $A$ (i.e., $A$ is a $G$-group), the $n$-th  $G$-commutator subgroup
$\Gamma^n_G(A) \subset A$ is the group generated by
$\{(ga^{-1})a\,|\, g\in G,\, a\in\Gamma_G^{n-1}(A)\}$, being
$\Gamma_G^0(A)=A$. The action is then nilpotent of nilpotency
order $r$, $\nil_GA=r$, if this is the smallest integer for which
$\Gamma^r_G(A)=\{1 \}$. The group $G$ also acts in each
$\Gamma_G^n(A)$ and
$\Gamma_G^m\bigl(\Gamma_G^n(A)\bigr)=\Gamma_G^{m+n}(A)$.

Statements of next sections shall heavily rely in the following
results:

\begin{lemma}\label{lema1}  Let $A$ be a $G$-group. Then:
\begin{itemize}
\item[(i)] $\Gamma_G^1(A)$ is a normal subgroup of $A$ and the $G$-action induced on $A/\Gamma_G^1(A)$ is trivial.
\item[(ii)] The quotient morphism  $A{\stackrel{q}{\longrightarrow}} A/\Gamma_G^1(A)$ is  equivariant and initial with
respect to trivial actions, i.e., every equivariant morphism
$A{\stackrel{f}{\longrightarrow}}H$, in which the $G$-action on
$H$ is trivial, factors uniquely through $q$.
\end{itemize}
\end{lemma}

\begin{proof} (i) is trivial. For (ii) observe that, for any $f$ as in the lemma, $\Gamma_G^1(A)\subset\ker f $.
\end{proof}

\begin{lemma}\label{lema2}  Let $A$ be a $G$-group. If $A$ is nilpotent then, for any $m$, $\Gamma_G^m(A)_\locp=\Gamma_G^m(A_\locp)$.
\end{lemma}

\begin{proof} Since $\Gamma_G^m(A)=\Gamma_G^{1}\bigl(\Gamma_G^{m-1}(A)\bigr)$, once we show that
$\Gamma_G^1(A)_\locp=\Gamma^1_G(A_\locp)$ an easy induction proves
the lemma. As localization  is an exact functor in the category of
nilpotent groups, the localization morphism $f:A\to A_{(p)}$
restricts to $f:\Gamma_G^1(A)\to \Gamma_G^1(A)_{(p)}$. Hence, we
may consider $\Gamma_G^1(A)_{(p)}$, as well as
$\Gamma_G^1(A_{(p)})$, as subgroups of $A_\locp$. Then, for any
$g\in G$ and $a\in A$, the trivial identity
$\bigl(gf(a)^{-1}\bigr)f(a)=f\bigl((ga^{-1})a\bigr)$ shows
equality of both subgroups.
\end{proof}

\begin{proposition}\label{uf} The group $G$ acts nilpotently on the nilpotent group $A$ if and only if $G$ acts
nilpotently on $A_\locp$ for $p$ any prime number or zero and all
these nilpotency orders are bounded.
\end{proposition}

\begin{proof} Assume $G$ acts nilpotently on $A$, i.e., $\Gamma_G^m(A)=\{1\}$ for some $m$. Hence, by Lemma \ref{lema2}
and for any $p$, $\Gamma_G^m(A_\locp)=\{1\}$.

Conversely, assume $\nil_GA_{(p)}\le m$, for all $p$ ($p$ a prime
number or $0$), and let $a$ be an element of $\Gamma_G^m(A)$. If
$a$ has finite order,  say it is a $q$-element, then it obviously
survives under the $q$-localization morphism
$\Gamma_G^m\bigl(A)\to\Gamma^m_G(A)_{(q)}$. For a general group,
elements of infinite order are not guaranteed to survive under
rationalization (for instance, the rationalization of the free
product of two finite groups is trivial while it contains elements
of infinite order). However for a nilpotent group, which is our
case, one can easily show by induction on the nilpotency order of
the group, that any element of infinite order is not sent to zero
under rationalization.  Taking into account, again by Lemma
\ref{lema2}, that $\Gamma^m_G(A)_{(p)}=\Gamma^m_G(A_{(p)})=\{1\}$,
it follows that $a=1$ and the proof is complete.
\end{proof}

\begin{proposition}\label{nuevolema} Let $G$ be a group acting on a finite nilpotent group $A$ in such a way that the
induced action on the Frattini factor $A/\Phi(A)$ is nilpotent.
Then, the $G$-action on $A$ is also nilpotent.
\end{proposition}

\begin{proof} Recall \cite[5.1]{Gorestein} that the Frattini subgroup of a group $A$, $\Phi(A)$, is defined to be the intersection of all its maximal proper subgroups. The Frattini factor of $A$ is $A/\Phi(A)$.
Observe in the first place that, since $\Phi(A)$ is a characteristic subgroup of $A$, i.e., it is invariant under any automorphism of $A$, $G$ in fact induces a
natural  action on the Frattini factor $A/\Phi(A)$ which, by
hypothesis, is nilpotent. Hence, since $A/\Phi(A)$ is nilpotent,
the induced action on
$\bigl(A/\Phi(A)\bigr)_\locp=A_\locp/\Phi(A)_\locp$ is also
nilpotent by Lemma \ref{lema2}. Next, observe that for any finite
group $A$, $\Phi(A)_\locp=\Phi(A_\locp)$. Indeed, this is
immediate from the definition taking into account that
localization commutes with limits,  in particular, with
intersections (see for instance \cite{hmr}). Therefore, we conclude that $G$ acts nilpotently on
$A_\locp/\Phi(A_\locp)$. Considering $\varphi\colon
G\to\aut\bigl(A_\locp/\Phi(A_\locp)\bigr)$ via this action, and
taking into account that $A_\locp/\Phi(A_\locp)$ is a finite
$p$-group, we may apply \cite[Corollary 5.3.3]{Gorestein} to
obtain that $\varphi(G)$ is also a $p$-group. But the action of a
$p$-group on another $p$-group is always nilpotent, and therefore
$G$ acts nilpotently on $A_\locp$. Since this is the case for any
$p$ and $A$ is finite we may apply  Proposition \ref{uf} and the
proposition follows.
\end{proof}

From these results we immediately deduce:

\bigskip
\noindent {\it Proof of Theorems \ref{dos} and \ref{uno}.} Apply
directly Propositions \ref{uf} and \ref{nuevolema} above to the
subgroup $G$ of $\eva (X)$ to obtain that $G$ acts nilpotently on
$\pi_{\le N}(X)$. Then, the result follows from Theorem
\ref{principal}. \hfill$\square$
\bigskip

Closely related to Proposition \ref{nuevolema}, we have the
following:

\begin{proposition}\label{propodos} Let $G$ be a group acting on an abelian $p$-group $A$ which has an exponent $p^n$. If
$G$ acts nilpotently on $A\otimes\bz/p$, then it does so on $A$
and $$ \nil_GA\le n\cdot\nil_GA\otimes\bz/p. $$
\end{proposition}
\begin{proof} Call $r=\nil_GA\otimes\bz/p$ and observe that $\Gamma_G^m(A\otimes\bz/p)=\Gamma_G^m(A)\otimes\bz/p$ for any
$m$. Therefore, since $\Gamma_G^r(A\otimes\bz/p) =0$,
$\Gamma_G^r(A)\subset pA$. Assume, as induction hypothesis, that
$\Gamma_G^{kr}(A)\subset p^kA$, for $k<n$. Hence, $$
\Gamma_G^{nr}(A)=\Gamma_G^{(n-1)r}\bigl(\Gamma_G^r(A)\bigr)\subset
\Gamma_G^{(n-1)r}(pA)= p\Gamma_G^{(n-1)r}(A)\subset p^nA. $$ Since
$A$ has $p^n$ as exponent, the proposition follows.
\end{proof}

\begin{proposition}\label{propouno} Let $A$ be a finite abelian
$p$-group and let $G\subset\aut(A)$ be such that $\sigma\otimes
{\bz/p}=1_{A\otimes\bz/p}$ for each $\sigma\in G$. Then $G$ is a
$p$-group.
\end{proposition}

\begin{proof}  As $A$ is a finite  abelian $p$-group,
the Frattini factor $A/\Phi(A)$ (respec. the projection
$A\rightarrow A/\Phi(A)$) is naturally  identified with
$A\otimes\bz/p$ (respec.\ the map $A\rightarrow A\otimes\bz/p$).
Now, if $G$ is  not a $p$-group, there exists a non trivial
$p^\prime$-automorphism $\sigma\in G$ which, by hypothesis and
using the identification above,  induces the identity on the
Frattini factor of $A$. But according to \cite[Theorem
5.1.4]{Gorestein}, the only $p^\prime$-automorphism that induces
the identity on the Frattini factor of a $p$-group is the
identity.  Thus $G$ must be a $p$-group.
\end{proof}

As an immediate consequence we get:

\begin{corollary}\label{corolario} In the conditions of the proposition above, the action of $G$ on $A$ is nilpotent.
\end{corollary}

\begin{proof} Indeed, recall that the action of a
$p$-group $H$ on another $p$-group is always nilpotent.
\end{proof}

In which follows $G$ is a group acting on another group $A$. Given
$g,h\in G$ and $a\in A$, we use the following usual notation: $$
[a,g]=a^{-1}(ga),\quad [g,a]=(ga^{-1})a,\quad
[g,h]=g^{-1}h^{-1}gh. $$ Hence, the following, which can be
considered as a variation of the Witt-Hall identity \cite[Theorem
5.1]{magnus}, is obtained by direct calculation.

\begin{lemma}\label{jo}
For any $f,g\in G$ and $b\in A$, the following identity holds:
$$\big[[f^{-1},g^{-1}],gb\big]b^{-1}\big[[g,b^{-1}], f\big]
b\big[[f,b],fgf^{-1}\big]=1.$$
\end{lemma}

\begin{lemma}\label{jo2} Let $H$ be a subgroup of $G$ and $K$ a normal subgroup of $H$. Then,
$$
\bigr[[H,K],A\bigr]\subset\langle\big[K,[H,A]\bigr],\bigl[H,[K,A]\bigr]\rangle.
$$
\end{lemma}

\begin{proof}
Making $f^{-1}=h$, $g^{-1}=k$ and $gb=a$ in Lemma \ref{jo},  it
follows that $$
\bigl[[h,k],a\bigr]=\bigl[h^{-1}k^{-1}h,[h^{-1},g^{-1}a]\bigr]
g^{-1}a^{-1}\bigl[h^{-1},[k^{-1},g^{-1}a]\bigr]g^{-1}a. $$ As $K$
is normal in $H$, $\bigl[h^{-1}k^{-1}h,[h^{-1},g^{-1}a]\bigr]\in
\big[K,[H,A]\bigr]$. On the other hand, as commutators are normal
subgroups,
$$g^{-1}a^{-1}\bigl[h^{-1},[k^{-1},g^{-1}a]\bigr]g^{-1}a\in
g^{-1}a^{-1} \bigl[H,[K,A]\bigr]g^{-1}a=\bigl[H,[K,A]\bigr], $$
and the lemma follows.
\end{proof}

\begin{lemma}\label{jo3} If the action of $G$ on $A$ is nilpotent, then
for each $n,m\ge 0$, $$ [\Gamma^n(G),\Gamma_G^m(A)]\subset
\Gamma_G^{n+m+1}(A). $$ In particular,
$
[\Gamma^n (G), A]\subset \Gamma_G^{n+1}(A).
$
\end{lemma}

\begin{proof} Set $\nil_GA=r$. If $m\ge r$ the assertion is obvious. Assume the lemma holds for all $n$ and $m\le 1$ and let us prove it for $m=0$ by
induction on $n$:

Trivially,
$
[\Gamma^0(G),\Gamma_G^0(A)]=[G,A]= \Gamma_G^{1}(A)$. Finally, $$
\begin{aligned}
{[\Gamma^n(G),A]}&=\bigl[[G,\Gamma^{n-1}(G)], A\bigr]\subset
\text{(By Lemma \ref{jo2})}\\ &\subset \langle
\bigl[\Gamma^{n-1}(G),[G,A]\bigr],
\bigl[G,[\Gamma^{n-1}(G),A]\bigr]\rangle\subset \text{(By
induction)}\\ &\subset\langle [\Gamma^{n-1}(G),\Gamma^1_G(A)],
[G,\Gamma^n_G(A)]\rangle\subset \text{(Again by induction)}\\
&\subset\Gamma_G^{n+1}(A).\\
\end{aligned}
$$
\end{proof}

\begin{proposition}\label{jo4}
Let $G$ be a subgroup of $\aut(A)$. Then, $\nil G\le \nil_GA-1$.
\end{proposition}

\begin{proof}
Assume $\nil_GA=r$. By Lemma \ref{jo3},
$[\Gamma^{r-1}(G),A]\subset \Gamma_G^r(A)=\{1\}$, and therefore
$\Gamma^{r-1}(G)=\{1\}$.
\end{proof}

\section{Self homotopy equivalence of spaces with local coefficients}

As stated in the Introduction, and following the notation and
approach of the standard reference \cite[Chap. VI.2]{white}, in
this section we consider self homotopy equivalences in the
homotopy category ${\cal L}^*$ of based spaces with local
coefficients. Observe that a self homotopy equivalence of an
object $(X,\calm)\in {\cal L}^*$ is given by  $(f,\Theta)\colon
(X,\calm)\to (X,\calm)$ in which $f\colon X\to X$ is a based
homotopy equivalence and $\Theta\colon \calm\to\calm $ is an
isomorphism of the coefficient system $\calm$. Note that such a
self equivalence $(f,\Theta)$ acts  in $\pi_*(X)$ by $\pi_*f$, in
$\calm$ by $\Theta$, and in $H^*(X;\calm)$ by $H^*(f,\Theta)$.

It is also convenient to  recall how cohomology classes with local
coefficients are represented by maps into the ``twisted
Eilenberg-MacLane space'' (see \cite[Chapter 5.2]{baues},
\cite{Gitler}, \cite{moller2} or  \cite{siegel} for precise
details). Let $K({\cal M},n)$ be a fixed realization of the
Eilenberg-MacLane space of type $(M_{x_0},n)$ being $M_{x_0}$ the
group of the system $\calm$ at the base point. On the other hand,
denote by $L({\cal M},n)$ the space obtained by applying the Borel
construction  to the universal fibration $\pi_1(X)\to\widetilde
K\stackrel{q}{\to} K(\pi_1(X),1)$ and the space $K({\cal M},n)$,
i.e.,

$$ L({\cal M},n)=\widetilde K\times_{\pi_1(X)} K({\cal M},n), $$
and it fits into the fibration $$ K(\calm,n)\longrightarrow
L(\calm,n)\stackrel{p}{\longrightarrow}K(\pi_1(X),1),\quad
p(a,b)=q(a). $$ Then, for a given space $X$, $H^n(X;{\cal M})$ is
in one to one correspondence with the set $[X,L({\cal
M},n)]_{K(\pi_1(X),1)}$ of homotopy classes of maps over
$K(\pi_1(X),1)$ from $X$ to $L({\cal M},n)$.

\bigskip

\noindent {\it Proof of Theorem \ref{local}.} To avoid excessive
notation we shall not distinguish between a homotopy class and a
map which represents it. First, observe that, if $(f,\Theta)$ is a
self homotopy equivalence of $(X,\calm)$ and $\alpha\colon X\to
L(\calm,n)$ represents a class of $H^*(X;\calm)$,
$H^n(f,\Theta)(\alpha)$ is represented by the map $$
X\stackrel{f}{\longrightarrow}
X\stackrel{\alpha}{\longrightarrow}L(\calm,n)\stackrel{\xi}{\longrightarrow}L(\calm,n)
$$ in which $\xi$ is defined by the action of $\Theta_{x_0}$ on
$M_{x_0}$. Explicitly, for  $(a,b)\in L(\calm,n)$,
$\xi(a,b)=(a,\overline{\Theta}_{x_0}b)$ where and
$\overline{\Theta}_{x_0}$ is the realization of $\Theta_{x_0}$.
Observe that $\xi$ is well defined as $\Theta\colon \calm\to\calm$
is a morphism of local coefficient systems.

Moreover, if $\alpha,\beta\colon X\to L(\calm,n)$ are in
$H^n(X;\calm)$, they coincides after composing with the fibration
$p\colon L(\calm,n)\to K(\pi(X),1)$, and therefore, for each $x\in
X$, $\alpha(x),\beta(x)$ live in the same fiber $K(\calm,n)$ of
$p$. Hence, $\alpha$ and $\beta$ can be added up on $K(\calm,n)$
and the resulting map $\alpha+\beta$ represents precisely their
sum as cohomology classes with twisted coefficients.

We shall prove the theorem by induction on the length of the
Postnikov decomposition of $X$. Assume  $X=K(\pi,m)$ and let
$(f,\Theta)\colon (X,\calm)\to(X,\calm)$ be a self equivalence.
Then, in view of the above, for any $n$-cohomology class
$\alpha\colon K(\pi,m)\to L(\calm,n)$,
$H^n(f,\Theta)(\alpha)-\alpha$ is represented by the map
$K(\pi,m)\to L(\calm,n)$ which, fiberwise, is
$\overline\Theta_{x_0}\alpha f-\alpha$. Writing
$\overline\Theta_{x_0}\alpha f-\alpha= \overline\Theta_{x_0}\alpha
f- \overline\Theta_{x_0}\alpha + \overline\Theta_{x_0}\alpha
-\alpha$ it is straightforward, using the nilpotency hypothesis,
to show that the $s$-th commutator of the action of $G$ on
$H^*(X;\calm)$ vanishes as long  as
$s\le\max\{\nil_G\pi_*(X),\nil_G\calm\}$.

Assume the theorem holds for $X=X^{r-1}$ and let $X=X^r$ be a
$r$-dimensional Postnikov piece.  Consider the Serre spectral
sequence with local coefficients on ${\cal M}$ associated to the
fibration $$ K(\pi_r(X),r)\to X\to X^{r-1}. $$ whose $E_2$-term is
$$ E_2^{*,*}=H^*\bigl(X^{r-1};{\cal H}^*(K(\pi_r(X),r);{\cal
M})\bigr). $$ Note that $G$ acts naturally in the base, total
space and fiber of this fibration, and hence, it does so in all
the terms of the spectral sequence.  The same argument used for
$r=1$ shows that  $G$ acts nilpotently on the local coefficient
system ${\cal H}^*(K(\pi_r(X),r);{\cal M})$ and therefore, by
induction hypothesis, $G$ acts nilpotently on
$H^*\bigl(X^{r-1};{\cal H}^*(K(\pi_r(X),r);{\cal M})\bigr)$.

As the spectral sequence converges, the action of $G$ on the
associated graded module of  $H^*(X;\calm)$ is nilpotent. Finally,
reasoning by induction on the filtration degree we deduce that the
$G$-action on  $H^*(X;\calm)$ is also nilpotent.\hfill$\square$

In particular, for any space $X$ and any $j$ we may consider the
local coefficient system given by $\pi_jX$. In this case, any self
homotopy equivalence $f\in\eva(X)$ can be seen as a self homotopy
equivalence $(f,\pi_jf^{-1})\in\eva(X;\pi_jX)$. Hence, any
subgroup of $\eva(X)$ may be considered as a subgroup of
$\eva(X,\pi_jX)$ which then acts naturally on $H^*(X;\pi_jX)$ when
considering local coefficients. In this context,  the theorem
above reads:

\begin{corollary}\label{noveas} Let $X$ be a finite Postnikov piece and  let $G$ be a subgroup of
$\eva(X)$ which acts nilpotently on $\pi_*(X)$. Then, for any $j$,
$G$ acts nilpotently on $H^*(X;\pi_j)$.\hfill$\square$
\end{corollary} This result is used in the proof of Theorem \ref{principal}  that we now present:

\bigskip
\noindent {\it Proof of Theorem \ref{principal}.} Consider the
restriction to $G$ of the exact sequence $
{1}\to\es(X)\to\eva(X)\to  \Pi_{i\le N}\text{aut}\,\pi_i(X) $: $$
{1}\to\es(X)\cap G\to G\to \Pi_{i\le N}\text{aut}\,\pi_i(X). $$
The image  of $G$ under this morphism, call it $\widetilde G$, is
a subgroup of automorphism of the group $\pi_{\le N}(X)$ in which
$G$ acts nilpotently by hypothesis. Then, by Proposition
\ref{jo4}, $\widetilde G$ is itself nilpotent and $\nil\widetilde
G< \nil_G\pi_{\le N}(X)$. Therefore, if we prove that $G$ acts
nilpotently  on $\es(X)$, then (see for instance \cite[Proposition
4.1]{hmr}) $G$ would be nilpotent and
\begin{equation}\label{nose1}
\nil G< \nil_G\es(X)+\nil_G\pi_{\le N}(X).
\end{equation} For that, observe in the first place that $[X,X]\cong[X^N,X^N]$, where
$X^N$ denotes the $N$-th Postnikov stage of $X$, and this
bijection restricts to an isomorphism $\es(X)\cong\es(X^N)$. On
the other hand, consider the exact sequence
\begin{equation}\label{nose2} {1}\to A_j\to \es(X^j)\to \es(X^{j-1})
\end{equation} where $\es(X^j)\to
\es(X^{j-1})$ is just the obvious restriction and $A_j$ its
kernel. Since $G$ acts on any $\es(X^j)$ and $\es(X^1)={1}$, it
will be enough to show that $G$ acts nilpotently on every $A_j$ to
conclude, by an easy induction, that it does so on
$\es(X^N)=\es(X)$.

By classical obstruction theory of liftings (see \cite[Chapter
6.6]{white}) recall that, for $j\ge 2$, there is a bijection
$\varphi\colon B_j\to H^j(X^j;\pi_j)$ where

\begin{itemize}
\item
The cohomology is taken with local coefficients.

\item $B_j$ is the set of homotopy classes of $[X^j,X^j]$ which restrict to  the identity on $X^{j-1}$, i.e., homotopy
classes of liftings of $X^j\to X^{j-1}$ to $X^j$.

\item $\varphi(g)=\delta(g,1)$ is the difference cochain of degree $j$ between $g$ and the identity on $X^j$.
   \end{itemize}
Recall also that, in general,
$\delta(g,f)=\delta(g,1)+\delta(1,f)$ and that
$\delta(gh,fh)=H^j(h)\bigl(\delta(g,f)\bigr)$. Moreover, if
$h\in\es(X^j)$, $\delta(hg,hf)$ is the image of $\delta(g,f)$
under the map $H^j(X^j;\pi_j)\to H^j(X^j;\pi_j)$ induced by $h$ on
$\pi_j$.

From now on, as in Corollary \ref{noveas}, any $f\in\eva(X)$, and
thus in $G$, shall be considered as a self homotopy class
$(f,\pi_jf^{-1})\in\eva(X,\pi_j)$. Hence, restricting $\varphi$ to
$A_j$ we obtain a map $\varphi\colon A_j\hookrightarrow
H^j(X^j;\pi_j)$ which is a $G$-map with respect to the action
$g\cdot f=g^{-1}fg$, $g\in G$, $f\in A_j$,  and the usual action
on $H^j(X^j;\pi_j)$: if $\alpha\in H^j(X_j,\pi_j)$, and $g\in G$,
$g\cdot\alpha$ is the cohomology class represented by the map

$$ X^j{\stackrel{g}{\longrightarrow}}
X^j{\stackrel{\alpha}{\longrightarrow}}
L(\pi_j,j){\stackrel{\xi}{\longrightarrow}} L(\pi_j,j),\qquad
\alpha\in H^j(X^j;\pi_j), $$ with $\xi\colon L(\pi_j,j)\to
L(\pi_j,j)$ induced by $\pi_*(g^{-1})$.

Moreover, this restriction is a group morphism. Indeed, given
$f,h\in A_j$,
$\varphi(fh)=\delta(fh,1)=\delta(f^{-1}fh,f^{-1})=\delta(h,f^{-1})=\delta(h,1)+\delta(1,f^{-1})=\delta(h,1)+\delta(f,1)=
\varphi(f)+\varphi(h)$. As an immediate consequence we then obtain
that $A_j$ is an (abelian!) subgroup of $H^j(X^j;\pi_j)$.

Finally, by Corollary \ref{noveas}, $G$ acts nilpotently on
$H^j(X^j;\pi_j)$ for any $j$. Hence, it does so on $A_j$ and $$
\nil_G\,A_j\le \nil_G\,H^j(X^j;\pi_j). $$ Thus, by induction on
$j$ using repeatedly \cite[Proposition 4.1]{hmr}, one easily sees
in view of (\ref{nose2}) that
$\nil_G\,\es(X)=\nil_G\,\es(X^N)\le\sum_{j=2}^N\nil_G\,A_j\le\sum_{j=2}^N\nil_G\,H^j(X^j;\pi_j)$
and the theorem follows.\hfill$\square$

\section{Groups which fix the homotopy groups}
In this section we establish Theorems \ref{tres} and \ref{cuatro}.

\bigskip
\noindent {\it Proof of Theorem \ref{tres}.}   Let
$\alpha\in\esp(X)$. Then, for each $i\le\dim X$, the morphism
$\pi_i(\alpha;\bz/p)\colon\pi_i(X,\bz/p)\to\pi_i(X,\bz/p)$ is just
the identity. On the other hand, the Universal Coefficients
Theorem for homotopy yields the following split short exact
sequence $$
0\longrightarrow\text{Ext}(\bz/p,\pi_{i+1}X)\longrightarrow\pi_i(X;\bz/p)\longrightarrow
\hom(\bz/p;\pi_iX)\longrightarrow 0. $$ Thus, both
$\text{Ext}(\bz/p,\pi_{i+1}\alpha)$ and $\hom(\bz/p;\pi_i\alpha)$
are the identity. But observe that
$\text{Ext}(\bz/p,\pi_{i+1}X)=\pi_{i+1}X\otimes\bz/p$ so that
$\pi_{i+1}\alpha\otimes\bz/p=1_{\pi_{i+1}X}$.

Hence, by Proposition \ref{propouno}, for each $i\le\dim X$, the
image of $\esp(X)$ in $\aut(\pi_iX)$ is a $p$-group and then, by
Corollary \ref{corolario}, the action of $\esp(X)$ on $\pi_iX$ is
nilpotent. Thus, by Theorem \ref{principal}, $\esp(X)$ is
nilpotent. On the other hand, notice that $\es(X)$ is precisely
the kernel of the obvious map $\esp(X)\to\Pi_{i\le\dim
X}\text{aut}\pi_i(X)$. Hence, as we just proved that the image of
this map is a $p$-group, $\esp(X)/\es(X)$ is a finite $p$-group,
and the proof is complete. \hfill$\square$

\bigskip

\noindent {\it Proof of Theorem \ref{cuatro}.} Write
$\pi_iX=\bz^{n_i}\oplus\bigl(\oplus_{p\,\,\text{prime}}T_p(\pi_iX)\bigr)$
in which $T_p(\pi_iX)$ is the group of $p$-torsion elements in
$\pi_iX$. Now, if $\alpha\in\cap_{p\,\,\text{prime}}\esp(X)$, then
for each $i\le\dim X$,
$\pi_i(\alpha)|_{T_p(\pi_iX)}\in\aut\bigl(T_p(\pi_iX)\bigr)$.  Let
$z_1,\ldots,z_{n_i}$ be generators of $\bz^{n_i}\subset\pi_iX$.
Then $\pi_i(\alpha)(z_k)=\sum_{j=1}^1{n_i}a_{k,j}z_j+\omega$,
where $\omega$ is the torsion part. But this element has to
coincide with $z_k$ mod $p$, for all prime $p$. Therefore, also
for any $p$, $ a_{k,j}=0(\text{mod}\,p) $ for $k\not=j$ and
$a_{k,k}=1(\text{mod}\, p)$, for $1\le k\le n_i$. The only
possible solution is $a_{k,j}=0$, $k\not=j$, and $a_{k,k}=1$. In
other words, $\pi_i(\alpha)(z_k)=z_k+\omega$, in which $\omega$ is
a torsion element.

Adding up, for any element $\gamma\in\pi_iX$,
$\pi_i(\alpha)\gamma-\gamma$ is a torsion element in $\pi_iX$.
This is equivalent to say that the $1$-commutators of the action
of $\cap_{p\,\,\text{prime}}\esp(X)$ on $\pi_iX$ live in the
torsion part of $\pi_iX$. However, by Corollary \ref{corolario},
the action of $\cap_{p\,\,\text{prime}}\esp(X)$ on the torsion
part is nilpotent and therefore,  the action on $\pi_iX$ is also
nilpotent. Apply Theorem \ref{principal} and the proof is
complete.\hfill$\square$

\begin{remark} {\em We end up by noting that the hypothesis of Theorem \ref{tres} are necessary. Indeed, consider
$X=K\bigl((\bz/2)^2,n\bigr)$ and observe that, for a prime $p$
different from $2$, $\esp(X)=\eva(X)=GL_2(\bz/2)\cong \Sigma_3$
which is not nilpotent.

On the other hand, take $X=K(\bz^2,n)$ for which
$\eva(X)=GL_2(\bz)$. In this case $\es(X)=\{1\}$ and, for any
prime $p$, $\esp(X)$ fits in the following short exact sequence $$
\{1\}\to\esp(X)\to GL_2(\bz)\to GL_2(\bz/p)\to\{1\} $$ where the
surjection $GL_2(\bz)\to GL_2(\bz/p)$ is just the mod-$p$
reduction. Hence $\esp(X)=\esp(X)/\es(X)$ is an infinite, non
nilpotent group.}
\end{remark}

\nocite{*}
\bibliographystyle{plain}

\end{document}